\documentclass[reqno,12pt]{article}
\usepackage{amsmath,amsfonts,amssymb,colordvi}
\def\Z{\mathbb{Z}}

\def\qed{\nopagebreak\hfill{\rule{4pt}{7pt}}}
\def\proof{ \noindent {\it{Proof.} \hskip 2pt}}

\newtheorem{theo}{Theorem}[section]

\newtheorem{lemm}[theo]{Lemma}

\makeatletter \@addtoreset{figure}{section} \makeatother
\makeatletter
\long\def\@makecaption#1#2{%
   \vskip 10\p@
   \setbox\@tempboxa\hbox{{#1}\ \ #2}%
   \ifdim \wd\@tempboxa >\hsize
       {#1}\ \ #2\par
   \else
       \hbox to\hsize{\hfil\box\@tempboxa\hfil}%
   \fi}
\makeatother
\newdimen\Squaresize \Squaresize=14pt
\newdimen\Thickness \Thickness=0.7pt

\def\Square#1{\hbox{\vrule width \Thickness
   \vbox to \Squaresize{\hrule height \Thickness\vss
      \hbox to \Squaresize{\hss#1\hss}
   \vss\hrule height\Thickness}
\unskip\vrule width \Thickness} \kern-\Thickness}

\def\Vsquare#1{\vbox{\Square{$#1$}}\kern-\Thickness}

\def\young#1{
\vbox{\smallskip\offinterlineskip \halign{&\Vsquare{##}\cr #1}}}

\def\moins{\raise 1pt\hbox{{$\scriptstyle -$}}}

\begin{document}
\begin{center}
{\large \bf Stable Equivalences of Giambelli Type Matrices of
Schur Functions}
\end{center}

\begin{center}
William Y. C. Chen$^{1}$ and Arthur
L. B. Yang$^{2}$\\[6pt]
Center for Combinatorics, LPMC\\
Nankai University, Tianjin 300071, P. R. China\\
Email: $^{1}${\tt chen@nankai.edu.cn}, $^{2}${\tt
yang@nankai.edu.cn}
\end{center}

\vspace{0.3cm} \noindent{\bf Abstract.}  By using cutting strips
and transformations on outside decompositions of a skew diagram,
we show that the Giambelli type matrices of a skew Schur function
are stably equivalent to each other over  symmetric functions. As
a consequence, the Jacobi-Trudi matrix and the dual Jacobi-Trudi
matrix are stably equivalent over symmetric functions. This gives
an affirmative answer to  an open problem posed by Kuperberg.



\section{Introduction}

Let $R$ be a commutative ring with unit. Recall that two matrices
$M$ and $M'$ over $R$ are called to be \emph{stably equivalent} to
each other, if and only if $M$ and $M'$ can be transformed from
each other by the following three fundamental matrix operations:
\begin{itemize}
\item[(i)] general row operation: $M\rightarrowtail AM=M'$;

\item[(ii)] general column operation: $M\rightarrowtail MA=M'$;

\item[(iii)] stabilization: $M\rightarrowtail
\begin{pmatrix}1 & 0\\ 0 & M\end{pmatrix}=M'$, and its inverse,
\end{itemize}
where $A$ is an invertible matrix over $R$.

This paper is motivated by Kuperberg's  problem \cite{KupEjc} on
the stable equivalence property between the Jacobi-Trudi matrix
and the dual Jacobi-Trudi matrix of skew Schur functions over the
ring $\Lambda$ of symmetric functions. We assume that the reader
is familiar with the notation and terminology on symmetric
functions in \cite{Stanley}. Given a partition $\lambda$ with
weakly decreasing components, let $\ell(\lambda)$ denote the
length of $\lambda$. The Jacobi-Trudi matrix for the skew Schur
function $s_{\lambda/\mu}$ is given by
\begin{equation}
J_{\lambda/\mu}=
\left(h_{\lambda_i-\mu_j-i+j}\right)_{i,j=1}^{\ell(\lambda)},
\end{equation}
where $h_k$ denotes the $k$-th complete symmetric function,
$h_0=1$ and $h_k=0$ for $k<0$. The dual Jacobi-Trudi matrix for
$s_{\lambda/\mu}$ is given by
\begin{equation}
D_{\lambda/\mu}=
\left(e_{\lambda'_i-\mu'_j-i+j}\right)_{i,j=1}^{\ell(\lambda')},
\end{equation}
where $\lambda'$ is the partition conjugate to $\lambda$, $e_k$
denotes the $k$-th elementary symmetric function, $e_0=1$ and
$e_k=0$ for $k<0$.

Kuperberg showed that the Jacobi-Trudi matrix and the dual
Jacobi-Trudi matrix are stably equivalent over the polynomial ring
\cite[Theorem 14]{KupEjc}.  He raised the question whether they
are  stably equivalent over the ring of symmetric functions. We
give an affirmative answer to this problem.

This paper is organized as follows. First we review some concepts
of outside decompositions for a given skew diagram. Utilizing the
cutting strips for a given edgewise connected skew shape as
introduced by Chen, Yan and Yang  \cite{CYY}, we demonstrate  how
a twist transformation changes the set of contents of the initial
boxes of border strips in an outside decomposition, and how it
changes the set of the contents of the terminal boxes. In Section
3,  we construct the canonical form of the Giambelli type matrix
of the skew Schur function assuming that the outside decomposition
is fixed. Using this canonical form we establish the stable
equivalence property of the Giambelli type matrix for the edgewise
connected skew diagram. In Section 4, we show that for a general
skew diagram the Jacobi-Trudi matrix and its dual form of Schur
functions are stably equivalent over the ring of symmetric
functions.

\section{Twist transformations}

Let $\lambda$ be a \emph{partition} of $n$ with $k$ parts, i.e.,
$\lambda=(\lambda_1,\lambda_2,\ldots,\lambda_k)$ where
$\lambda_1\geq\lambda_2\geq\ldots\geq\lambda_k>0$ and
$\lambda_1+\lambda_2+\ldots+\lambda_k=n$. A \emph{Young diagram}
of $\lambda$ may be defined as the set of points $(i,j)\in \Z^2$
such that $1\leq j\leq\lambda_i$ and $1\leq i\leq k$. A Young
diagram can also be represented in the plane by an array of boxes
justified from the top and left corner with $k$ rows and
$\lambda_i$ boxes in row $i$. A box $(i,j)$ in the diagram is the
box in row $i$ from the top and column $j$ from the left. The
content of $(i,j)$, denoted $\tau((i,j))$, is given by  $j-i$.
Given two partitions $\lambda$ and $\mu$, we say that
$\mu\subseteq\lambda$ if $\mu_i\leq \lambda_i$ for all $i$. If
$\mu\subseteq\lambda$, we define a \emph{skew partition}
$\lambda/\mu$, whose Young diagram is obtained from the Young
diagram of $\lambda$ by peeling off the Young diagram of $\mu$
from the upper left corner. The \emph{conjugate} of a skew
partition $\lambda/\mu$, which we denote by $\lambda'/\mu'$, is
defined to be the transpose of the skew diagram $\lambda/\mu$.

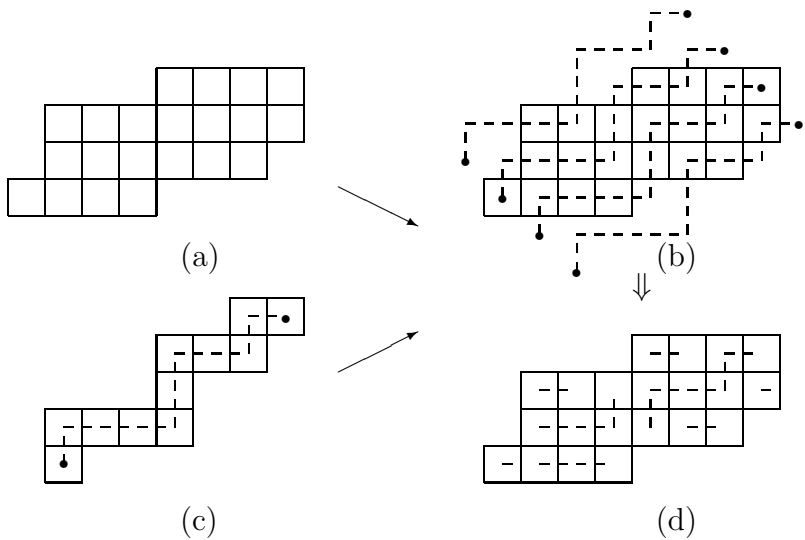
\begin{figure}[h,t]
\begin{center}
\setlength{\unitlength}{20pt}
\begin{picture}(16,9.6)

\put(4,0){(c)}  \put(13,0){(d)} \put(4,5){(a)}  \put(13,5){(b)}

 \put(0.75,5.95){\line(1,0){2.8}}
\put(0.75,6.65){\line(1,0){4.9}} \put(1.45,7.35){\line(1,0){4.9}}
\put(1.45,8.05){\line(1,0){4.9}} \put(3.55,8.75){\line(1,0){2.8}}

\put(0.75,5.95){\line(0,1){0.7}} \put(1.45,5.95){\line(0,1){2.1}}
\put(2.15,5.95){\line(0,1){2.1}} \put(2.85,5.95){\line(0,1){2.1}}
\put(3.55,5.95){\line(0,1){2.8}} \put(4.25,6.65){\line(0,1){2.1}}
\put(4.95,6.65){\line(0,1){2.1}} \put(5.65,6.65){\line(0,1){2.1}}
\put(6.35,7.35){\line(0,1){1.4}}

\put(1.45,0.9){\line(1,0){0.7}} \put(1.45,1.6){\line(1,0){2.8}}
\put(1.45,2.3){\line(1,0){2.8}} \put(3.55,3.0){\line(1,0){2.1}}
\put(3.55,3.7){\line(1,0){2.8}} \put(4.95,4.4){\line(1,0){1.4}}

\put(1.45,0.9){\line(0,1){1.4}} \put(2.15,0.9){\line(0,1){1.4}}
\put(2.85,1.6){\line(0,1){0.7}} \put(3.55,1.6){\line(0,1){2.1}}
\put(4.25,1.6){\line(0,1){2.1}} \put(4.95,3.0){\line(0,1){1.4}}
\put(5.65,3.0){\line(0,1){1.4}} \put(6.35,3.7){\line(0,1){0.7}}

\multiput(1.8,1.95)(0.35,0){6}{\line(1,0){0.175}}
\multiput(3.9,3.35)(0.35,0){4}{\line(1,0){0.175}}
\multiput(5.3,4.05)(0.35,0){2}{\line(1,0){0.175}}

\multiput(1.8,1.25)(0,0.35){2}{\line(0,1){0.175}}
\multiput(3.9,1.95)(0,0.35){4}{\line(0,1){0.175}}
\multiput(5.3,3.35)(0,0.35){2}{\line(0,1){0.175}}

\put(1.8,1.25){\circle*{0.125}} \put(6,4){\circle*{0.125}}

\put(7,6.5){\vector(2,-1){1.5}} \put(7,3){\vector(2,1){1.5}}
\put(12.55,4.45){$\Downarrow$}

\put(9.75,5.95){\line(1,0){2.8}} \put(9.75,6.65){\line(1,0){4.9}}
\put(10.45,7.35){\line(1,0){4.9}}
\put(10.45,8.05){\line(1,0){4.9}}
\put(12.55,8.75){\line(1,0){2.8}}

\put(9.75,5.95){\line(0,1){0.7}} \put(10.45,5.95){\line(0,1){2.1}}
\put(11.15,5.95){\line(0,1){2.1}}
\put(11.85,5.95){\line(0,1){2.1}}
\put(12.55,5.95){\line(0,1){2.8}}
\put(13.25,6.65){\line(0,1){2.1}}
\put(13.95,6.65){\line(0,1){2.1}}
\put(14.65,6.65){\line(0,1){2.1}}
\put(15.35,7.35){\line(0,1){1.4}}

\multiput(0.35,4.65)(0.7,-0.7){4}{\multiput(9.05,3.05)(0.35,0){6}{\line(1,0){0.175}}}
\multiput(0.35,4.65)(0.7,-0.7){4}{\multiput(11.15,4.45)(0.35,0){4}{\line(1,0){0.175}}}
\multiput(0.35,4.65)(0.7,-0.7){4}{\multiput(12.55,5.15)(0.35,0){2}{\line(1,0){0.175}}}

\multiput(0.35,4.65)(0.7,-0.7){4}{\multiput(9.05,2.35)(0,0.35){2}{\line(0,1){0.175}}}
\multiput(0.35,4.65)(0.7,-0.7){4}{\multiput(11.15,3.05)(0,0.35){4}{\line(0,1){0.175}}}
\multiput(0.35,4.65)(0.7,-0.7){4}{\multiput(12.55,4.45)(0,0.35){2}{\line(0,1){0.175}}}

\multiput(9.4,6.975)(0.7,-0.7){4}{\circle*{0.125}}
\multiput(13.6,9.775)(0.7,-0.7){4}{\circle*{0.125}}

\put(9.75,0.9){\line(1,0){2.8}} \put(9.75,1.6){\line(1,0){4.9}}
\put(10.45,2.3){\line(1,0){4.9}} \put(10.45,3){\line(1,0){4.9}}
\put(12.55,3.7){\line(1,0){2.8}}

\put(9.75,0.9){\line(0,1){0.7}} \put(10.45,0.9){\line(0,1){2.1}}
\put(11.15,0.9){\line(0,1){2.1}} \put(11.85,0.9){\line(0,1){2.1}}
\put(12.55,0.9){\line(0,1){2.8}} \put(13.25,1.6){\line(0,1){2.1}}
\put(13.95,1.6){\line(0,1){2.1}} \put(14.65,1.6){\line(0,1){2.1}}
\put(15.35,2.3){\line(0,1){1.4}}

\put(10.1,1.25){\line(1,0){0.175}}
\multiput(10.8,1.25)(0.35,0){4}{\line(1,0){0.175}}
\multiput(10.8,1.95)(0.35,0){4}{\line(1,0){0.175}}
\multiput(10.8,2.65)(0.35,0){2}{\line(1,0){0.175}}
\multiput(12.9,3.35)(0.35,0){2}{\line(1,0){0.175}}
\multiput(13.6,1.95)(0.35,0){2}{\line(1,0){0.175}}
\multiput(12.9,2.65)(0.35,0){4}{\line(1,0){0.175}}
\multiput(14.3,3.35)(0.35,0){2}{\line(1,0){0.175}}

\multiput(12.2,1.95)(0,0.35){2}{\line(0,1){0.175}}
\multiput(12.9,1.95)(0,0.35){2}{\line(0,1){0.175}}
\multiput(14.3,2.65)(0,0.35){2}{\line(0,1){0.175}}
\put(15,2.65){\line(1,0){0.175}}

\end{picture}
\end{center}
\caption{The cutting strip of an outside decomposition}
\label{fig-2}
\end{figure}

A skew diagram $\lambda/\mu$ is \emph{connected} if it can be
regarded as a union of an edgewise connected set of boxes, where
two boxes are said to be edgewise connected if they share a common
edge. A \emph{border strip} is a connected skew diagram with no
$2\times 2$ block of boxes. If no two boxes lie in the same row,
we call such a border strip  a \emph{vertical border strip}. If no
two boxes lie in the same column, we call such a border strip  a
\emph{horizontal border strip}. An \emph{outside decomposition} of
$\lambda/\mu$ is a partition of the boxes of $\lambda/\mu$ into
pairwise disjoint border strips such that every border strip in
the decomposition has a starting box on the left or bottom
perimeter of the diagram and an ending box on the right or top
perimeter of the diagram, see Figure \ref{fig-2} (d). This concept
was used by Hamel and Goulden \cite{HG} to give a lattice path
proof for the Giambelli type determinant formulas of the skew
Schur function.

Recall that a \emph{diagonal} with content $c$ of $\lambda/\mu$ is
the set of all the boxes in $\lambda/\mu$  having content $c$.
Starting from the lower left corner of the skew diagram
$\lambda/\mu$, we use  consecutive integers $1,2,\ldots,d$ to
number these diagonals. Chen, Yan and Yang \cite{CYY} obtained the
following characterization of outside decompositions of a skew
shape.

\begin{theo}[{\cite[Theorem 2.2]{CYY}}]
Suppose that $\lambda/\mu$ is an edgewise connected skew partition
and has $d$ diagonals. Then there is a one-to-one correspondence
between the outside decompositions of $\lambda/\mu$ and border
strips with $d$ boxes.
\end{theo}

For each outside decomposition $\Pi$, the corresponding border
strip $T$ is called the \emph{cutting strip} of $\Pi$ in
\cite{CYY}, which is given by the rule: for $i=1, 2, \ldots, d-1$,
the relative position between the $i$-th box and the $(i+1)$-th
box in $T$ coincides with the relative position between the two
boxes in the same border strip of $\Pi$, one of which is on the
$i$-th diagonal and the other on the $(i+1)$-th diagonal, see
Figure \ref{fig-2}.

Notice that the relative position between the $i$-th box and the
$(i+1)$-th box of the border strip imposes an up or right
direction to the $i$-th box according to the $(i+1)$-th box lies
above or to the right of the $i$-th box. Throughout this paper, we
will read the diagram from the bottom right corner to the top left
corner. In the same manner, each diagonal may be endowed with a
direction for a given outside decomposition. From the cutting
strip characterization of outside decompositions, one can obtain
any outside decomposition from another by a sequence of basic
transformations of changing the directions of the boxes on a
diagonal, which corresponds to the operation of changing the
direction of a box in the cutting strip. This transformation is
called {\em the twist transformation} on border strips.

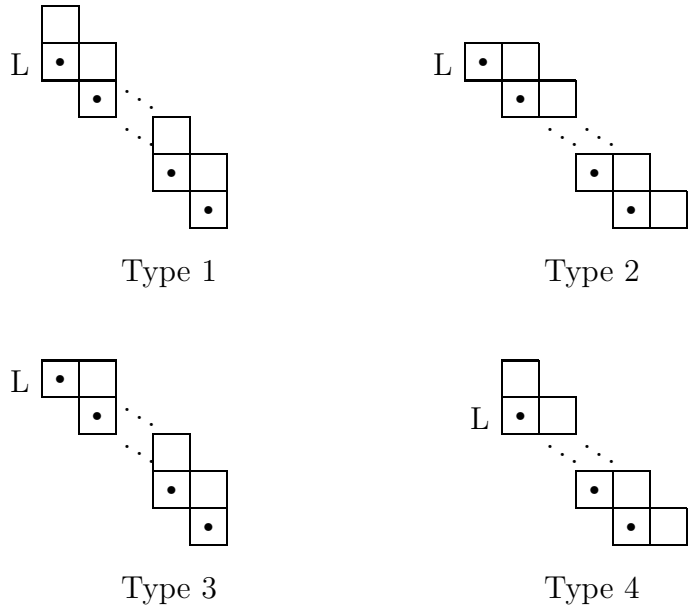
\begin{figure}[h,t]
\setlength{\unitlength}{20pt}
\begin{center}
\begin{picture}(15,12)

\put(3.8,7){\line(1,0){0.7}} \put(3.1,7.7){\line(1,0){1.4}}
\put(3.1,8.4){\line(1,0){1.4}} \put(3.1,9.1){\line(1,0){0.7}}
\put(1.7,9.1){\line(1,0){0.7}} \put(1,9.8){\line(1,0){1.4}}
\put(1,10.5){\line(1,0){1.4}} \put(1,11.2){\line(1,0){0.7}}

\put(0.4,9.9){L} \put(8.4,9.9){L}

\put(1,9.8){\line(0,1){1.4}} \put(1.7,9.1){\line(0,1){2.1}}
\put(2.4,9.1){\line(0,1){1.4}} \put(3.1,7.7){\line(0,1){1.4}}
\put(3.8,7){\line(0,1){2.1}} \put(4.5,7){\line(0,1){1.4}}

\put(2.5,9.2){$\ddots$}\put(2.5,8.5){$\ddots$}

\put(1.35,10.15){\circle*{0.125}} \put(4.15,7.35){\circle*{0.125}}
\put(3.45,8.05){\circle*{0.125}} \put(2.05,9.45){\circle*{0.125}}

\put(11.8,7){\line(1,0){1.4}} \put(11.1,7.7){\line(1,0){2.1}}
\put(11.1,8.4){\line(1,0){1.4}} \put(9.7,9.1){\line(1,0){1.4}}
\put(9,9.8){\line(1,0){2.1}} \put(9,10.5){\line(1,0){1.4}}

\put(9,9.8){\line(0,1){0.7}} \put(9.7,9.1){\line(0,1){1.4}}
\put(10.4,9.1){\line(0,1){1.4}} \put(11.1,7.7){\line(0,1){0.7}}
\put(11.1,9.1){\line(0,1){0.7}} \put(11.8,7){\line(0,1){1.4}}
\put(12.5,7){\line(0,1){1.4}} \put(13.2,7){\line(0,1){0.7}}

\put(11.2,8.5){$\ddots$} \put(10.5,8.5){$\ddots$}

\put(9.35,10.15){\circle*{0.125}}
\put(12.15,7.35){\circle*{0.125}}
\put(11.45,8.05){\circle*{0.125}}
\put(10.05,9.45){\circle*{0.125}}

\put(2.5,6){Type 1 } \put(10.5,6){Type 2}

\put(0.4,3.9){L} \put(9.1,3.2){L}

\put(3.8,1){\line(1,0){0.7}} \put(3.1,1.7){\line(1,0){1.4}}
\put(3.1,2.4){\line(1,0){1.4}} \put(3.1,3.1){\line(1,0){0.7}}
\put(1.7,3.1){\line(1,0){0.7}} \put(1,3.8){\line(1,0){1.4}}
\put(1,4.5){\line(1,0){1.4}}

\put(1,3.8){\line(0,1){0.7}} \put(1.7,3.1){\line(0,1){1.4}}
\put(2.4,3.1){\line(0,1){1.4}} \put(3.1,1.7){\line(0,1){1.4}}
\put(3.8,1){\line(0,1){2.1}} \put(4.5,1){\line(0,1){1.4}}

\put(2.5,3.2){$\ddots$}\put(2.5,2.5){$\ddots$}

\put(1.35,4.15){\circle*{0.125}} \put(4.15,1.35){\circle*{0.125}}
\put(3.45,2.05){\circle*{0.125}} \put(2.05,3.45){\circle*{0.125}}

\put(11.8,1){\line(1,0){1.4}} \put(11.1,1.7){\line(1,0){2.1}}
\put(11.1,2.4){\line(1,0){1.4}} \put(9.7,3.1){\line(1,0){1.4}}
\put(9.7,3.8){\line(1,0){1.4}} \put(9.7,4.5){\line(1,0){0.7}}

\put(9.7,3.1){\line(0,1){1.4}} \put(10.4,3.1){\line(0,1){1.4}}
\put(11.1,1.7){\line(0,1){0.7}} \put(11.1,3.1){\line(0,1){0.7}}
\put(11.8,1){\line(0,1){1.4}} \put(12.5,1){\line(0,1){1.4}}
\put(13.2,1){\line(0,1){0.7}}

\put(11.2,2.5){$\ddots$} \put(10.5,2.5){$\ddots$}

\put(12.15,1.35){\circle*{0.125}}
\put(11.45,2.05){\circle*{0.125}}
\put(10.05,3.45){\circle*{0.125}}

\put(2.5,0){Type 3 } \put(10.5,0){Type 4}
\end{picture}
\end{center}
\caption{Four possible types of diagonals of $\lambda/\mu$}
\label{4types}
\end{figure}

Let $\lambda/\mu$ be an edgewise connected skew shape. Let $L$ be
the diagonal of $\lambda/\mu$ consisting of the boxes with content
$i$. Note that $L$ must be one of the four possible diagonal types
classified by whether the first diagonal box has a box at the top,
and whether the last diagonal box has a box on the right. These
four types are depicted by Figure \ref{4types}.

Given an outside decomposition $\Pi=(\theta_1, \theta_2, \ldots,
\theta_m)$ of $\lambda/\mu$ and a strip $\theta$ in $\Pi$, we
denote the content of the initial box of $\theta$ and the content
of the terminal box of $\theta$ respectively by $p(\theta)$ and
$q(\theta)$. Let $\phi$ be the cutting strip of $\Pi$. It is known
that $\theta$ can be regarded as the segment of $\phi$ with the
initial content $p(\theta)$ and the terminal content $q(\theta)$
\cite{CYY}, denoted $\phi[p(\theta), q(\theta)]$.

Given two skew diagrams $I$ and $J$, let $I\blacktriangleright J$
be the diagram obtained by gluing the lower left-hand corner box
of diagram $J$ to the right of the upper right-hand corner box of
diagram $I$, and let $I\uparrow J$ be the diagram obtained by
gluing the lower left-hand corner box of diagram $J$ to the top of
the upper right-hand corner box of diagram $I$. Suppose that the
strip $\theta$ has a box of $L$, then $\theta$ can be written as
$\phi[p(\theta), i]\blacktriangleright \phi[i+1,q(\theta)]$ if $L$
has the right direction, and $\theta$ can be written as
$\phi[p(\theta), i]\uparrow \phi[i+1,q(\theta)]$ if $L$ has the up
direction.

Let $\omega_i$ denote the twist transformation acting on an
outside decomposition $\Pi$ by changing the direction of the
diagonal $L$. Let
\begin{eqnarray}
{\rm Init}(\Pi) & = & \{p(\theta_1),p(\theta_2),\ldots,p(\theta_m)\},\\
{\rm Term}(\Pi) & = &
\{q(\theta_1),q(\theta_2),\ldots,q(\theta_m)\}.
\end{eqnarray} The following theorem describes the actions of $\omega_i$
on ${\rm Init}(\Pi)$ and ${\rm Term}(\Pi)$.

\begin{theo} \label{mainthm} Given an outside decomposition $\Pi$, let $\Pi'$ be the
outside decomposition obtained from $\Pi$ by applying the twist
transformation $\omega_i$. Then we have
\begin{itemize}

\item[(a)] $i\not\in {\rm Term}(\Pi), i+1\not\in {\rm Init}(\Pi),
{\rm Init}(\Pi')={\rm Init}(\Pi)\cup\{i+1\}$ and ${\rm
Term}(\Pi')={\rm Term}(\Pi)\cup\{i\}$, or

\item[(b)] $i\in {\rm Term}(\Pi), i+1\in {\rm Init}(\Pi), {\rm
Init}(\Pi')={\rm Init}(\Pi)\setminus\{i+1\}$ and ${\rm
Term}(\Pi')={\rm Term}(\Pi)\setminus\{i\}$, or

\item[(c)] $i\in {\rm Term}(\Pi), i+1\not\in {\rm Init}(\Pi), {\rm
Init}(\Pi')={\rm Init}(\Pi)$ and ${\rm Term}(\Pi')={\rm
Term}(\Pi)$, or

\item[(d)] $i\not\in {\rm Term}(\Pi), i+1\in {\rm Init}(\Pi), {\rm
Init}(\Pi')={\rm Init}(\Pi)$ and ${\rm Term}(\Pi')={\rm
Term}(\Pi)$.

\end{itemize}
\end{theo}
\proof Suppose that $L$ has $r$ boxes. Since the twist
transformation $\omega_i$ only changes the strips which contain a
box in $L$, we may suppose that $\theta_{i_t}, 1\leq t\leq r$, is
the strip in $\Pi$ that contains the $t$-th diagonal box in $L$.
Without loss of generality we may assume that the diagonal boxes
have the up direction, since we can reverse the transformation
process for the case when the diagonal boxes have the right
direction.

Let $\phi'$ be the cutting strip corresponding to the outside
decomposition $\Pi'$. Now we see the changes of ${\rm Init}(\Pi)$
and ${\rm Term}(\Pi)$ under the action of the twist transformation
$\omega_i$ according to the type of $L$:

\begin{itemize}

\item[(a)] If $L$ is of Type 1, then  we have $i\not\in {\rm
Term}(\Pi)$ and $i+1\not\in {\rm Init}(\Pi)$. As illustrated in
Figure \ref{D1}, under the operation of $\omega_i$, the strip
\[ \theta_{i_1} =
\phi[p(\theta_{i_1}),q(\theta_{i_1})]=\phi[p(\theta_{i_1}),i]\uparrow\phi[i+1,q(\theta_{i_1})]\]
 breaks into two strips
 \[
\phi'[p(\theta_{i_1}),q(\theta_{i_2})]=\phi[p(\theta_{i_1}),i]\blacktriangleright\phi[i+1,q(\theta_{i_2})]
\quad \mbox{ and} \quad \phi'[i+1,q(\theta_{i_1})].\]
 If $r > 1$ then the last
strip
\[ \theta_{i_r} =
\phi[p(\theta_{i_r}),q(\theta_{i_r})]=\phi[p(\theta_{i_r}),i]\uparrow\phi[i+1,q(\theta_{i_r})]\]
will be cut off into $\phi'[p(\theta_{i_r}),i]$, and the other
strips
\[ \theta_{i_t} =
\phi[p(\theta_{i_t}),q(\theta_{i_t})]=\phi[p(\theta_{i_t}),i]\uparrow\phi[i+1,q(\theta_{i_t})],\
2 \leq t \leq r-1, \] will be twisted into
\[\phi'[p(\theta_{i_t}),q(\theta_{i_{t+1}})]=\phi[p(\theta_{i_t}),i]\blacktriangleright\phi[i+1,q(\theta_{i_{t+1}})].\]
Thus
$$
 {\rm
Init}(\Pi')={\rm Init}(\Pi)\cup\{i+1\} \mbox { and } {\rm
Term}(\Pi')={\rm Term}(\Pi)\cup\{i\}.$$

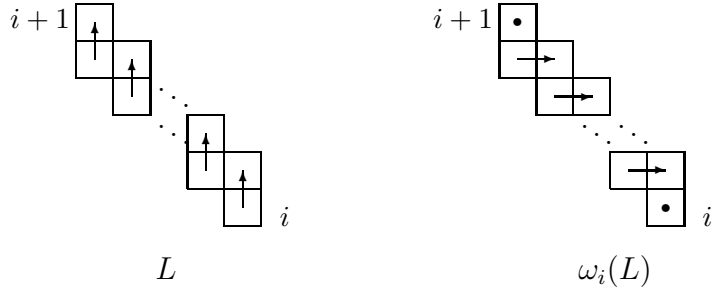
\begin{figure}[h,t]
\begin{center}
\setlength{\unitlength}{20pt}
\begin{picture}(14,5.5)

\put(3.8,1){\line(1,0){0.7}}\put(3.1,1.7){\line(1,0){1.4}}
\put(3.1,2.4){\line(1,0){1.4}}\put(3.1,3.1){\line(1,0){0.7}}
\put(1.7,3.1){\line(1,0){0.7}}\put(1,3.8){\line(1,0){1.4}}
\put(1,4.5){\line(1,0){1.4}}\put(1,5.2){\line(1,0){0.7}}

\put(1,3.8){\line(0,1){1.4}}\put(1.7,3.1){\line(0,1){2.1}}
\put(2.4,3.1){\line(0,1){1.4}}\put(3.1,1.7){\line(0,1){1.4}}
\put(3.8,1){\line(0,1){2.1}}\put(4.5,1){\line(0,1){1.4}}

\put(4.15,1.35){\vector(0,1){0.7}}\put(3.45,2.05){\vector(0,1){0.7}}
\put(2.05,3.45){\vector(0,1){0.7}}\put(1.35,4.15){\vector(0,1){0.7}}

\put(4.85,1.0){\small $i$} \put(-0.25,4.75){\small $i+1$}

\put(2.5,3.2){$\ddots$}\put(2.5,2.5){$\ddots$}

\put(11.8,1){\line(1,0){0.7}}\put(11.1,1.7){\line(1,0){1.4}}
\put(11.1,2.4){\line(1,0){1.4}}
\put(9.7,3.1){\line(1,0){1.4}}\put(9,3.8){\line(1,0){2.1}}
\put(9,4.5){\line(1,0){1.4}}\put(9,5.2){\line(1,0){0.7}}

\put(9,3.8){\line(0,1){1.4}}\put(9.7,3.1){\line(0,1){2.1}}
\put(10.4,3.1){\line(0,1){1.4}}\put(11.1,1.7){\line(0,1){0.7}}
\put(11.1,3.1){\line(0,1){0.7}}
\put(11.8,1){\line(0,1){1.4}}\put(12.5,1){\line(0,1){1.4}}

\put(11.45,2.05){\vector(1,0){0.7}}
\put(10.05,3.45){\vector(1,0){0.7}}\put(9.35,4.15){\vector(1,0){0.7}}

\put(11.2,2.5){$\ddots$}\put(10.5,2.5){$\ddots$}

\put(9.35,4.85){\circle*{0.125}} \put(12.15,1.35){\circle*{0.125}}

\put(12.85,1.0){\small $i$} \put(7.75,4.75){\small $i+1$}

\put(2.5,0){$L$ } \put(10.5,0){$\omega_i(L)$}
\end{picture}
\end{center}
\caption{$\omega_i$ acts on a Type 1 diagonal $L$ }\label{D1}
\end{figure}

\item[(b)] If $L$ is of Type 2, then we have $i\in {\rm
Term}(\Pi)$ and $i+1\in {\rm Init}(\Pi)$. Let $\theta_{i_{r+1}}$
be the unique strip of $\Pi$ with the initial content $i+1$. Under
the operation of $\omega_i$, the strip $\theta_{i_1}
=\phi[p(\theta_{i_1}),i]$ becomes a part of the new strip
$$\phi'[p(\theta_{i_1}),q(\theta_{i_2})].$$ The strip
$\theta_{i_{r+1}}=\phi[i+1,q(\theta_{i_{r+1}})]$ becomes a part of
the new strip
$$\phi'[p(\theta_{i_{r}}),q(\theta_{i_{r+1}})]=\phi[p(\theta_{i_r}),i]\blacktriangleright\phi[i+1,q(\theta_{i_{r+1}})].$$ For $2\leq
t\leq r-1$, the strips
\[ \theta_{i_t} =
\phi[p(\theta_{i_t}),q(\theta_{i_t})]=\phi[p(\theta_{i_t}),i]\uparrow\phi[i+1,q(\theta_{i_t})]
\] will be twisted into
\[\phi'[p(\theta_{i_t}),q(\theta_{i_{t+1}})]=\phi[p(\theta_{i_t}),i]\blacktriangleright\phi[i+1,q(\theta_{i_{t+1}})].\]
Thus
$$
 {\rm
Init}(\Pi')={\rm Init}(\Pi)\setminus\{i+1\} \mbox { and } {\rm
Term}(\Pi')={\rm Term}(\Pi)\setminus\{i\}.$$

\item[(c)] If $L$ is of Type 3, then we have $i\in {\rm
Term}(\Pi)$ and $i+1\not\in {\rm Init}(\Pi)$. Under the operation
of $\omega_i$, the first strip
$\theta_{i_1}=\phi[p(\theta_{i_1}),i]$ becomes
\[\phi'[p(\theta_{i_1}),q(\theta_{i_2})]=\phi[p(\theta_{i_1}),i]\blacktriangleright\phi[i+1,q(\theta_{i_2})].\]
If $r>1$, the last strip
\[\theta_{i_r}=\phi[p(\theta_{i_r}),q(\theta_{i_r})]=\phi[p(\theta_{i_r}),i]\uparrow\phi[i+1,q(\theta_{i_r})]\]
will be cut off into $\phi'[p(\theta_{i_r}),i]$, and the other
strips
\[\theta_{i_t}=\phi[p(\theta_{i_t}),q(\theta_{i_t})]=\phi[p(\theta_{i_t}),i]\uparrow\phi[i+1,q(\theta_{i_t})],\]
will be twisted into
\[\phi'[p(\theta_{i_t}),q(\theta_{i_{t+1}})]=\phi[p(\theta_{i_t}),i]\blacktriangleright\phi[i+1,q(\theta_{i_{t+1}})],\ 2\leq t\leq r-1.\]
Thus
$${\rm Init}(\Pi')={\rm Init}(\Pi) \mbox{ and } {\rm Term}(\Pi')={\rm
Term}(\Pi).$$

\item[(d)] If $L$ is of Type 4, then we have $i\not\in {\rm
Term}(\Pi)$ and $i+1\in {\rm Init}(\Pi)$. Let $\theta_{i_{r+1}}$
be the unique strip of $\Pi$ with the initial content $i+1$. Under
the operation $\omega_i$, the first strip
$$\theta_{i_1}=\phi[p(\theta_{i_1}),q(\theta_{i_1})]=\phi[p(\theta_{i_1}),i] \uparrow \phi[i+1,q(\theta_{i_1})]$$ breaks into
two strips
\[\phi'[p(\theta_{i_1}),q(\theta_{i_2})]=\phi[p(\theta_{i_1}),i] \blacktriangleright \phi[i+1,q(\theta_{i_2})]
\mbox{ and } \phi'[i+1,q(\theta_{i_1})].\] The strip
$\theta_{i_{r+1}}$ becomes a part of the new strip
\[\phi'[p(\theta_r),q_{\theta_{r+1}}]=\phi[p(\theta_r),r]\blacktriangleright \phi[i+1,q(\theta_{i_{r+1}})].\]
The other strips
\[ \theta_{i_t} =
\phi[p(\theta_{i_t}),q(\theta_{i_t})]=\phi[p(\theta_{i_t}),i]\uparrow\phi[i+1,q(\theta_{i_t})],\
2 \leq t \leq r-1, \] will be twisted into
\[\phi'[p(\theta_{i_t}),q(\theta_{i_{t+1}})]=\phi[p(\theta_{i_t}),i]\blacktriangleright\phi[i+1,q(\theta_{i_{t+1}})].\]
Thus $${\rm Init}(\Pi')={\rm Init}(\Pi) \mbox{ and } {\rm
Term}(\Pi')={\rm Term}(\Pi).$$

\end{itemize}

\qed

\section{Giambelli type matrices for connected shapes}

By using the lattice path methodology, Hamel and Goulden \cite{HG}
give a combinatorial proof for the Giambelli type determinant
formulas of the skew Schur function. In this section, we prove the
stable equivalence of the Giambelli type matrices of the Schur
function  indexed by an edgewise connected skew partition
$\lambda/\mu$.

Given an outside decomposition $\Pi=(\theta_1, \theta_2, \ldots,
\theta_m)$ of $\lambda/\mu$ and a strip $\theta$ in $\Pi$, let
$\phi$ be the cutting strip of $\Pi$. Recall that the strip
$\theta$ coincides with the segment $\phi[p(\theta),q(\theta)]$ of
$\phi$. Following the treatment of \cite{CYY}, given any two
contents $p,q$ we may define the strip $\phi[p,q]$  as follows:
\begin{enumerate}

\item[(i)] If $p \leq q$, then $\phi[p,q]$ is the segment of
$\phi$ starting with the box having content $p$ and ending with
the box having content $q$;

\item[(ii)] If $p = q+1$, then $\phi[p,q]$ is the empty strip
$\emptyset$.

\item[(iii)] If $p > q+1$, then $\phi[p,q]$ is undefined.

\end{enumerate}

Hamel and Goulden proved the following result.

\begin{theo}[{\cite[Theorem 3.1]{HG}}]\label{schur-det}
The skew Schur function $s_{\lambda/\mu}$ can be evaluated by the
following determinant:
\begin{equation}\label{eq1}
D(\Pi)= \det(s_{\phi[p(\theta_j),q(\theta_i)]})_{i,j=1}^{m}
\end{equation}
where $s_{\emptyset}=1$ and $s_{undefined}=0$.
\end{theo}

Let us denote the Giambelli type matrix   in \eqref{eq1} by
$M(\Pi)$. Chen, Yan and Yang \cite{CYY} have obtained the
canonical form $C(\Pi)=(s_{\phi[p_i,q_j]})_{i,j=1}^{m}$ of
$M(\Pi)$, where the sequence $(p_1, p_2, \ldots, p_m)$ is the
decreasing reordering of $(p(\theta_1), p(\theta_2), \ldots,
p(\theta_m))$ and $(q_1, q_2, \ldots, q_m)$ is the decreasing
reordering of $(q(\theta_1), q(\theta_2), \ldots, q(\theta_m))$.
It is clear that if $s_{[p_i,q_j]}=0$ then $s_{[p_i,\
q_{j\,'}]}=0$ and $s_{[p_{i\,'},\ q_j]}=0$ for $j\leq j\,'\leq m$
and $1\leq i\,'\leq i$.

Since $M(\Pi)$ and $C(\Pi)$ can be obtained from each other by
permutations of rows and columns. Thus we have

\begin{lemm}\label{key-lem}
For an outside decomposition $\Pi$ of the skew diagram
$\lambda/\mu$, the two matrices $M(\Pi)$ and $C(\Pi)$ are stably
equivalent over the ring $\Lambda$ of symmetric functions.
\end{lemm}

In order to show that the two Giambelli type matrices $M(\Pi)$ and
$M(\Pi')$ are stably equivalent over $\Lambda$, it suffices to
prove that their canonical forms $C(\Pi)$ and $C(\Pi')$ are stably
equivalent. To this end, we need the following lemma:

\begin{lemm}[\cite{LP2, Sch, Z}] \label{keylem1}
Let $I$ and $J$ be two skew diagrams. Then
\begin{equation}\label{eqv-rl}
s_{I}s_{J}=s_{I \blacktriangleright J}+s_{I \uparrow J}.
\end{equation}
\end{lemm}

We now come to the main theorem of this paper:
\begin{theo}\label{main-thm}
Let $\Pi$ and $\Pi'$ be two outside decompositions of the edgewise
connected skew diagram $\lambda/\mu$. Then the Giambelli type
matrices $M(\Pi)$ and $M(\Pi')$ are stably equivalent over the
ring $\Lambda$ of symmetric functions.
\end{theo}

\proof By Lemma \ref{key-lem}, we only need to prove that $C(\Pi)$
and $C(\Pi')$ are stably equivalent over $\Lambda$. Since any two
outside decompositions can be transformed from each other by a
sequence of twist transformations, it suffices to prove the case
when $\Pi'=\omega_i(\Pi)$ for any twist transformation $\omega_i$.
Let $\phi$ be the cutting strip of $\Pi$, and let $\phi'$ be the
cutting strip of $\Pi'$. Without loss of generality, we assume
that the box of content $i$ in $\phi$ has the up direction. Thus
the box of content $i$ in $\phi'$ has the right direction. By
Theorem \ref{mainthm}, we only need to consider the stable
equivalence between $C(\Pi)$ and $C(\Pi')$. Here are four cases:

\begin{enumerate}

\item[(a)] $i\not\in {\rm Term}(\Pi), i+1\not\in {\rm Init}(\Pi),
{\rm Init}(\Pi')={\rm Init}(\Pi)\cup\{i+1\}$ and ${\rm
Term}(\Pi')={\rm Term}(\Pi)\cup\{i\}$. Suppose that $k$ and $k'$
are the two indices such that
$$p_k>i+1 \mbox{ and } p_{k+1}<i+1; \mbox{ while } q_{k'}>i \mbox{ and } q_{k'+1}<i.$$
Then the canonical matrix $C(\Pi)$ has the following form
$$
\begin{pmatrix} s_{\phi[p_1,q_1]} & \cdots &
s_{\phi[p_1,q_{k'}]} & 0 &
\cdots & 0\\
\vdots & \vdots & \vdots & \vdots &
\vdots & \vdots\\
s_{\phi[p_k,q_1]} & \cdots & s_{\phi[p_k,q_{k'}]} & 0 & \cdots &
0\\
s_{\phi[p_{k+1},i]\uparrow \phi[i+1, q_1]} & \cdots &
s_{\phi[p_{k+1},i]\uparrow\phi[i+1, q_{k'}]} &
s_{\phi[p_{k+1},q_{k'+1}]} &
\cdots & s_{\phi[p_{k+1},q_{m}]}\\
\vdots & \vdots & \vdots & \vdots &
\vdots & \vdots\\
s_{\phi[p_m,i]\uparrow\phi[i+1,q_1]} & \cdots &
s_{\phi[p_m,i]\uparrow\phi[i+1,q_{k'}]} & s_{\phi[p_m,q_{k'+1}]} &
\cdots & s_{\phi[p_m,q_{m}]}
\end{pmatrix},
$$
and the canonical matrix $C(\Pi')$ has the following form
$$
\begin{pmatrix} s_{\phi[p_1,q_1]} & \cdots &
s_{\phi[p_1,q_{k'}]} & 0 & 0 &
\cdots & 0\\
\vdots & \vdots & \vdots & \vdots & \vdots &
\vdots & \vdots\\
s_{\phi[p_k,q_1]} & \cdots & s_{\phi[p_k,q_{k'}]} & 0 & 0 & \cdots
&
0\\
s_{\phi[i+1,q_1]} & \cdots & s_{\phi[i+1,q_{k'}]} & 1 & 0 & \cdots
&
0\\
s_{\phi[p_{k+1},i]\blacktriangleright \phi[i+1, q_1]} & \cdots &
s_{\phi[p_{k+1},i]\blacktriangleright \phi[i+1, q_{k'}]} &
s_{\phi[p_{k+1},i]} & s_{\phi[p_{k+1},q_{k'+1}]} &
\cdots & s_{\phi[p_{k+1},q_{m}]}\\
\vdots & \vdots & \vdots & \vdots &
\vdots & \vdots & \vdots\\
s_{\phi[p_m,i]\blacktriangleright \phi[i+1,q_1]} & \cdots &
s_{\phi[p_m,i]\blacktriangleright \phi[i+1,q_{k'}]} &
s_{\phi[p_{m},i]} & s_{\phi[p_m,q_{k'+1}]} & \cdots &
s_{\phi[p_m,q_{m}]}
\end{pmatrix}.
$$
For $j:1\leq j\leq k'$ subtracting the $(k'+1)$-th column of
$C(\Pi)$ multiplied by $s_{\phi[i+1, q_{j}]}$ from the $j$-th
column, then for $j:k+2\leq j\leq m+1$, subtracting the $(k+1)$-th
row multiplied by $s_{\phi[p_{j-1},i]}$ from the $j$-th row, we
get the following matrix due to Lemma \ref{keylem1}
$$
\begin{pmatrix} s_{\phi[p_1,q_1]} & \cdots &
s_{\phi[p_1,q_{k'}]} & 0 & 0 &
\cdots & 0\\
\vdots & \vdots & \vdots & \vdots & \vdots &
\vdots & \vdots\\
s_{\phi[p_k,q_1]} & \cdots & s_{\phi[p_k,q_{k'}]} & 0 & 0 & \cdots
&
0\\
0 & \cdots & 0 & 1 & 0 & \cdots &
0\\
-s_{\phi[p_{k+1},i]\uparrow \phi[i+1, q_1]} & \cdots &
-s_{\phi[p_{k+1},i]\uparrow \phi[i+1, q_{k'}]} & 0 &
s_{\phi[p_{k+1},q_{k'+1}]} &
\cdots & s_{\phi[p_{k+1},q_{m}]}\\
\vdots & \vdots & \vdots & \vdots &
\vdots & \vdots & \vdots\\
-s_{\phi[p_m,i]\uparrow \phi[i+1,q_1]} & \cdots &
-s_{\phi[p_m,i]\uparrow \phi[i+1,q_{k'}]} & 0 &
s_{\phi[p_m,q_{k'+1}]} & \cdots & s_{\phi[p_m,q_{m}]}
\end{pmatrix}.
$$
By multiplying $-1$ for the last $m-k$ rows and the last $m-k'$
columns, then permuting rows and columns, and the inverse
operation of stabilization, we find that the above matrix is
stably equivalent to $C(\Pi)$ over the ring $\Lambda$ of symmetric
functions. Thus $C(\Pi)$ and $C(\Pi')$ are stably equivalent over
$\Lambda$.

\item[(b)] $i\in {\rm Term}(\Pi), i+1\in {\rm Init}(\Pi), {\rm
Init}(\Pi')={\rm Init}(\Pi)\setminus\{i+1\}$ and ${\rm
Term}(\Pi')={\rm Term}(\Pi)\setminus\{i\}$. In this case, we only
need to reverse the process of the operations of case (a), where
$\omega_i$ is regarded as a transformation from the right
direction to the up direction. Notice that each inverse operation
is still over the ring $\Lambda$ of symmetric functions. Thus
$C(\Pi)$ and $C(\Pi')$ are stably equivalent over $\Lambda$.

\item[(c)] $i\in {\rm Term}(\Pi), i+1\not\in {\rm Init}(\Pi), {\rm
Init}(\Pi')={\rm Init}(\Pi)$ and ${\rm Term}(\Pi')={\rm
Term}(\Pi)$. Suppose that $k$ and $k'$ are the two indices such
that
$$p_k>i+1 \mbox{ and } p_{k+1}<i+1; \mbox{ while } q_{k'}=i.$$
Then the canonical matrix $C(\Pi)$ has the following form
$$
\begin{pmatrix} s_{\phi[p_1,q_1]} & \cdots &
s_{\phi[p_1,q_{k'-1}]} & 0 & 0 &
\cdots & 0\\
\vdots & \vdots & \vdots & \vdots &
\vdots & \vdots & \vdots\\
s_{\phi[p_k,q_1]} & \cdots & s_{\phi[p_k,q_{k'-1}]} & 0 & 0 &
\cdots &
0\\
s_{\phi[p_{k+1},i]\uparrow \phi[i+1, q_1]} & \cdots &
s_{\phi[p_{k+1},i]\uparrow\phi[i+1, q_{k'-1}]} &
s_{\phi[p_{k+1},i]} & s_{\phi[p_{k+1},q_{k'+1}]} &
\cdots & s_{\phi[p_{k+1},q_{m}]}\\
\vdots & \vdots & \vdots & \vdots & \vdots &
\vdots & \vdots\\
s_{\phi[p_m,i]\uparrow\phi[i+1,q_1]} & \cdots &
s_{\phi[p_m,i]\uparrow\phi[i+1,q_{k'-1}]} & s_{\phi[p_{m},i]} &
s_{\phi[p_m,q_{k'+1}]} & \cdots & s_{\phi[p_m,q_{m}]}
\end{pmatrix},
$$
and the canonical matrix $C(\Pi')$ has the following form
$$
\begin{pmatrix} s_{\phi[p_1,q_1]} & \cdots &
s_{\phi[p_1,q_{k'-1}]} & 0 & 0 &
\cdots & 0\\
\vdots & \vdots & \vdots & \vdots &
\vdots & \vdots & \vdots\\
s_{\phi[p_k,q_1]} & \cdots & s_{\phi[p_k,q_{k'-1}]} & 0 & 0 &
\cdots &
0\\
s_{\phi[p_{k+1},i]\blacktriangleright \phi[i+1, q_1]} & \cdots &
s_{\phi[p_{k+1},i]\blacktriangleright \phi[i+1, q_{k'-1}]} &
s_{\phi[p_{k+1},i]} & s_{\phi[p_{k+1},q_{k'+1}]} &
\cdots & s_{\phi[p_{k+1},q_{m}]}\\
\vdots & \vdots & \vdots & \vdots & \vdots &
\vdots & \vdots\\
s_{\phi[p_m,i]\blacktriangleright \phi[i+1,q_1]} & \cdots &
s_{\phi[p_m,i]\blacktriangleright \phi[i+1,q_{k'-1}]} &
s_{\phi[p_{m},i]} & s_{\phi[p_m,q_{k'+1}]} & \cdots &
s_{\phi[p_m,q_{m}]}
\end{pmatrix}.
$$
For $j:1\leq j\leq k'-1$ subtracting the $k'$-th column of
$C(\Pi)$ multiplied by $s_{\phi[i+1, q_{j}]}$ from the $j$-th
column, and then multiplying $-1$ for the last $m-k$ rows and the
last $m-k'+1$ columns, we obtain the matrix $C(\Pi)$. This implies
that $C(\Pi)$ and $C(\Pi')$ are stably equivalent over $\Lambda$.

\item[(d)] $i\not\in {\rm Term}(\Pi), i+1\in {\rm Init}(\Pi), {\rm
Init}(\Pi')={\rm Init}(\Pi)$ and ${\rm Term}(\Pi')={\rm
Term}(\Pi)$. We omit the proof here since it is similar to Case
(c).
\end{enumerate}
By summarizing, we have completed the proof. \qed

\section{Jacobi-Trudi matrices}

In this section we will prove that the Jabobi-Trudi matrix and the
dual Jacobi-Trudi matrix are stably equivalent over the ring
$\Lambda$ of symmetric functions for a general skew partition
$\lambda/\mu$. Theorem \ref{main-thm} states that this is true
when $\lambda/\mu$ is edgewise connected, where we do not allow
the existence of   empty strips in the outside decomposition
$\Pi$. The Jacobi-Trudi matrix $J_{\lambda/\mu}$ corresponds to
the Giambelli type matrix $M(\Pi)$  when the cutting strip $\phi$
of $\Pi$ is a horizontal border strip, and the dual Jacobi-Trudi
matrix $D_{\lambda/\mu}$ corresponds to the matrix $M(\Pi)$ when
$\phi$ is a vertical border strip.

For a general skew partition $\lambda/\mu$, we need to be more
careful when dealing with the empty strip. Let
$c_{min}=-\lambda_1'+1$ and $c_{max}=\lambda_1-1$. Let $\phi_h$
(or $\phi_e$) be the horizontal (resp. vertical) border strip
starting with the box having content $c_{min}$ and ending with the
box having content $c_{max}$. Let
$\Pi_h=(\theta_1,\cdots,\theta_{\ell(\lambda)})$ be the horizontal
outside decomposition of $\lambda/\mu$, where $\theta_i$ is a
horizontal strip of row $i$ from the $(\mu_i+1)$-th box to the
$\lambda_i$-th box. When $\lambda_i=\mu_i$, we take $\theta_i$ as
the empty strip. Clearly, each $\theta_i$ corresponds to a
substrip $\phi_h[\mu_i-i+1,\lambda_i-i]$ of $\phi_h$. Now we see
that the Jacobi-Trudi matrix $J_{\lambda/\mu}$ coincides with the
Giambelli type matrix $M(\Pi_h)$ defined in \eqref{eq1}.
Similarly, let $\Pi_e=(\theta_1',\cdots,\theta_{\lambda_1}')$ be
the vertical outside decomposition of $\lambda/\mu$, where
$\theta_i'$ is a vertical strip of column $i$ from the
$\lambda'_i$-th box to the $(\mu_i'+1)$-th box. When
$\lambda_i'=\mu_i'$, we take $\theta_i'$ as the empty strip.
Clearly, each $\theta_i'$ corresponds to a substrip
$\phi_e[-\lambda_i'+i,-\mu_i'+i-1]$ of $\phi_e$. Then the dual
Jacobi-Trudi matrix $D_{\lambda/\mu}$ coincides with the Giambelli
type matrix $M(\Pi_e)$. The following lemma is straightforward.

\begin{lemm}\label{trival-trans}
Let $\lambda/\mu$ be a partition with $\lambda_1'=\mu_1'$. Let
$\rho/\nu$ be the skew partition by removing the first column of
the skew diagram $\lambda/\mu$. Then the Jacobi-Trudi matrices of
$\lambda/\mu$ and $\rho/\nu$ are stably equivalent over $\Lambda$,
and so are the dual Jacobi-Trudi matrices.
\end{lemm}

Therefore, we may assume that $\lambda_1'\neq \mu_1'$. Let $\Pi$
be an outside decomposition of $\lambda/\mu$, and let $\phi$ be
the cutting strip of $\Pi$. For $i: c_{min}\leq i\leq c_{max}$,
let $\omega_i$ denote the twist transformation at the box of
content $i$ from the right direction to the up direction. Now we
define the outside decomposition $\omega_i(\Pi)$ by the following
rule:

\begin{enumerate}
\item[(a')] If $\lambda/\mu$ has both a box with content $i$ and a
box with content $i+1$, then let $\omega_i(\Pi)=\Pi\setminus
\Pi^{(i)}\cup \omega_i(\Pi^{(i)})$, where $\Pi^{(i)}$ is the
outside decomposition of the edgewise connected region of
$\lambda/\mu$ which has a box with content $i$ and
$\omega_i(\Pi^{(i)})$ is defined as in Section 2.

\item[(b')] If $\lambda/\mu$ has a box with content $i$ and but it
does not have a box with content $i+1$, then let
$\omega_i(\Pi)=\Pi$.

\item[(c')] If $\lambda/\mu$ neither have a box with content $i$
nor have a box with content $i+1$ while it has a box with content
less than $i$, then put $\omega_i(\Pi)=\Pi\cup \{\phi[i+1,i]\}$ if
$\phi[i+1,i]\not\in\Pi$, otherwise put $\omega_i(\Pi)=\Pi\setminus
\{\phi[i+1,i]\}$.

\item[(d')] If $\lambda/\mu$ has a box with content $i+1$ and a
box with content less than $i$, but it does not have a box with
content $i$, then let $\omega_i(\Pi)=\Pi$.

\end{enumerate}

The following lemma is a direct verification of the action of
$\omega_i$ on outside decompositions

\begin{lemm} \label{trans} Let $\lambda/\mu$ be a skew partition with $\lambda_1'\neq\mu_1'$.
Let $\Pi_h$ and $\Pi_e$ be the horizontal outside decomposition and the vertical
outside decomposition of $\lambda/\mu$ respectively. Then
\begin{equation}
\Pi_e=\omega_{c_{max}-1}(\omega_{c_{max}-2}(\cdots
(\omega_{c_{min}}({\Pi_h}))\cdots)).
\end{equation}
\end{lemm}

We now reach the following conclusion as an answer to Kuperberg's
problem \cite[Question 15]{KupEjc}.

\begin{theo}
For a skew partition $\lambda/\mu$, the Jacobi-Trudi matrix
$J_{\lambda/\mu}$ and $D_{\lambda/\mu}$ are stably equivalent over
the ring of symmetric functions.
\end{theo}

\noindent
 \proof Due to  Lemma \ref{key-lem}, we only need to
prove that the canonical matrices $C(\Pi_h)$ and $C(\Pi_e)$ are
stably equivalent over $\Lambda$. Due to Lemma \ref{trival-trans},
we only deal with the case of $\lambda_1'\neq\mu_1'$. By Lemma
\ref{trans}, it suffices to prove that $C(\Pi)$ and
$C(\omega_i(\Pi))$ are stably equivalent under any twist
transformation $\omega_i$ of the above four cases.

Let ${\rm Init(\Pi)}=\{p_1,p_2,\ldots,p_m\}$ and ${\rm
Term(\Pi)}=\{q_1,q_2,\ldots,q_m\}$ be strictly decreasing. Now we
see the transformations between the matrices according to the type
of $\omega_i$.

If $\omega_i$ is of type (c'), then the proof is similar to the
proof of case (a) and (b) in Theorem \ref{main-thm}.

For the case of $\omega_i$ being of type (a'), the stably
equivalent transformation will be one of the cases of the proof of
Theorem \ref{main-thm}.

If $\omega_i$ is of type (b'), then $i\in {\rm Term(\Pi)}$. Now
the proof is similar to the proof of case (c) in Theorem
\ref{main-thm}.

If $\omega_i$ is of type (d'), then $i+1\in {\rm Init(\Pi)}$. Now
the proof is similar to the proof of case (d) in Theorem
\ref{main-thm}.

Combining all the cases, we have completed the proof. \qed

\noindent {\textbf{Remark.}} The above proof only gives the stably
equivalent transformations from the Jacobi-Trudi matrix to the
dual Jacobi-Trudi matrix. In fact, we can also transform the dual
Jacobi-Trudi matrix into the Jacobi-Trudi matrix.

\begin{figure}[h,t]
\setlength{\unitlength}{1pt}
\begin{center}
\begin{picture}(186,85)
\put(40,85){\line(1,0){120}} \put(40,65){\line(1,0){120}}
\put(40,45){\line(1,0){100}} \put(40,25){\line(1,0){60}}
\put(40,5){\line(1,0){20}}

\put(40,5){\line(0,1){80}} \put(60,5){\line(0,1){80}}
\put(80,25){\line(0,1){60}} \put(100,25){\line(0,1){60}}
\put(120,45){\line(0,1){40}} \put(140,45){\line(0,1){40}}
\put(160,65){\line(0,1){20}}

\Gray{ \put(40.5,5.5){\rule{19pt}{19pt}}
\put(140.5,65.5){\rule{19pt}{19pt}}
\put(120.5,65.5){\rule{19pt}{19pt}}
\put(120.5,45.5){\rule{19pt}{19pt}}
 }
\end{picture}
\caption{The skew partition $(6,5,3,1)/(4,4,3)$}\label{fig-d}
\end{center}
\end{figure}
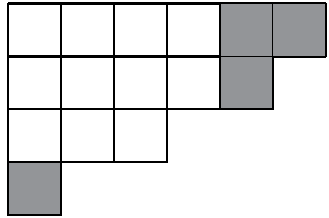

For instance, we take $\lambda/\mu=(6,5,3,1)/(4,4,3)$ to
illustrate the proof of the above theorem, see Appendix. The skew
diagram $\lambda/\mu=(6,5,3,1)/(4,4,3)$ is shown in Figure
\ref{fig-d}.

\vspace{.2cm} \noindent{\bf Acknowledgments.} This work was done
under the auspices of the 973 Project on Mathematical
Mechanization, the Ministry of Education, the Ministry of Science
and Technology, and the National Science Foundation of China.

\newpage
\section*{Appendix}

Note that the Jacobi-Trudi matrix is the transpose of the first
Giambelli type matrix, and the dual Jacobi-Trudi matrix is the
transpose of the last Giambelli type matrix. Here we use $[p,q]$
denote the corresponding border strip in the outside
decomposition. The dots in the matrix represent $0$.

\begin{table}[htb]
    \begin{center}
    \begin{tabular}{c|c}
         ${\mbox{Cutting strip and} \atop \mbox{outside decomposition}}$ & ${\mbox{Canonical form of}\atop \mbox{Giambelli type matrix}}$\\
         \hline
         $\begin{array}{c}
         \young{-3 & -2 & -1 & 0 &
        1 & 2 & 3 & 4 & 5\cr}\\
         \{[4,5],[3,3],[1,0],[-3,-3]\}
         \end{array}$
          & $\begin{pmatrix}
        s_2 & 1 & \cdot & \cdot \\
        s_3 & s_1 & \cdot & \cdot \\
        s_5 & s_3 & 1 & \cdot \\
        s_9 & s_7 & s_4 & s_1
        \end{pmatrix}$ \\
         \hline
    $\begin{array}{c}
         \young{-2 & -1 & 0 &
        1 & 2 & 3 & 4 & 5\cr
        {-3}\cr}\\
         \{[4,5],[3,3],[1,0],[-3,-3]\}
         \end{array}$
         & $
         \begin{array}{c}
         \\
         \begin{pmatrix}
        s_2 & 1 & \cdot & \cdot \\
        s_3 & s_1 & \cdot & \cdot \\
        s_5 & s_3 & 1 & \cdot \\
        s_{81} & s_{61} & s_{31} & s_1
        \end{pmatrix}
        \end{array}$
        \\
        \hline
$\begin{array}{c}
         \young{-1 & 0 &
        1 & 2 & 3 & 4 & 5\cr
        {-2}\cr
        -3\cr} \\
         \{[4,5],[3,3],[1,0],[-1,-2],[-3,-3]\}
         \end{array}$ &
        $\begin{pmatrix}
        s_2 & 1 & \cdot & \cdot &\cdot \\
        s_3 & s_1 & \cdot & \cdot & \cdot \\
        s_5 & s_3 & 1 & \cdot & \cdot \\
        s_{7} & s_5 & s_2 & 1 & \cdot\\
        s_{71^2} & s_{51^2} & s_{21^2} & s_{1^2} & s_1
        \end{pmatrix}$\\
        \hline
$\begin{array}{c}
         \young{0 &
        1 & 2 & 3 & 4 & 5 \cr
        {-1}\cr
        -2\cr
        -3\cr}\\
         \{[4,5],[3,3],[1,0],[0,-1], \\
         {[-1,-2]},[-3,-3]\}
         \end{array}$ &
        $\begin{pmatrix}
        s_2 & 1 & \cdot & \cdot & \cdot &\cdot \\
        s_3 & s_1 & \cdot & \cdot & \cdot & \cdot \\
        s_5 & s_3 & 1 & \cdot &  \cdot & \cdot \\
        s_6 & s_4 & s_1 & 1 &  \cdot & \cdot \\
        s_{61} & s_{41} & s_{1^2} & s_1 &  1 & \cdot\\
        s_{61^3} & s_{41^3} & s_{1^4} & s_{1^3} & s_{1^2} & s_1
        \end{pmatrix}$\\
        \hline
$\begin{array}{c}
         \young{
        1 & 2 & 3 & 4 & 5 \cr
        0 \cr
        {-1}\cr
        -2\cr
        -3\cr}\\
         \{[4,5],[3,3],[0,-1],
         {[-1,-2]},[-3,-3]\}
         \end{array}$ &
        $\begin{pmatrix}
        s_2 & 1 & \cdot & \cdot &\cdot \\
        s_3 & s_1  & \cdot & \cdot & \cdot \\
        s_{51} & s_{31}  & 1 &  \cdot & \cdot \\
        s_{51^2} & s_{31^2}  & s_1 &  1 & \cdot\\
        s_{51^4} & s_{31^4}  & s_{1^3} & s_{1^2} & s_1
        \end{pmatrix}$\\
        \hline
    \end{tabular}
    \end{center}
\end{table}

Continuing to the twist transformation, we have the following
\begin{table}[htb]
    \begin{center}
    \begin{tabular}{c|c}
         ${\mbox{Cutting strip and} \atop \mbox{outside decomposition}}$ & ${\mbox{Canonical form of}\atop \mbox{Giambelli type matrix}}$\\
        \hline
$\begin{array}{c}
         \young{
        2 & 3 & 4 & 5 \cr
        1 \cr
        0 \cr
        {-1}\cr
        -2\cr
        -3\cr}\\
         \{[4,5],[3,3],[2,1], \\
         {[0,-1]},
         {[-1,-2]},[-3,-3]\}
         \end{array}$ &
        $\begin{pmatrix}
        s_2 & 1 & \cdot & \cdot & \cdot &\cdot \\
        s_3 & s_1 & \cdot & \cdot & \cdot & \cdot \\
        s_4 & s_2 & 1 & \cdot & \cdot & \cdot\\
        s_{41^2} & s_{21^2} & s_{1^2}  & 1 &  \cdot & \cdot \\
        s_{41^3} & s_{21^3}  & s_{1^3} & s_1 &  1 & \cdot\\
        s_{41^5} & s_{21^5}  & s_{1^5} & s_{1^3} & s_{1^2} & s_1
        \end{pmatrix}$\\
        \hline
$\begin{array}{c}
         \young{
        3 & 4 & 5 \cr
        2\cr
        1 \cr
        0 \cr
        {-1}\cr
        -2\cr
        -3\cr}\\
         \{[4,5],[3,3],[2,1], \\
         {[0,-1]},
         {[-1,-2]},[-3,-3]\}
         \end{array}$ &
        $\begin{pmatrix}
        s_2 & 1 & \cdot & \cdot & \cdot &\cdot \\
        s_3 & s_1 & \cdot & \cdot & \cdot & \cdot \\
        s_{31} & s_{1^2} & 1 & \cdot & \cdot & \cdot\\
        s_{31^3} & s_{1^4} & s_{1^2}  & 1 &  \cdot & \cdot \\
        s_{31^4} & s_{1^5}  & s_{1^3} & s_1 &  1 & \cdot\\
        s_{31^6} & s_{1^7}  & s_{1^5} & s_{1^3} & s_{1^2} & s_1
        \end{pmatrix}$\\
        \hline
$\begin{array}{c}
         \young{
        4 & 5 \cr
        3\cr
        2\cr
        1 \cr
        0 \cr
        {-1}\cr
        -2\cr
        -3\cr}\\
         \{[3,5],[2,1],
         {[0,-1]},
         {[-1,-2]},[-3,-3]\}
         \end{array}$ &
        $\begin{pmatrix}
        s_{21} & \cdot & \cdot & \cdot & \cdot \\
        s_{21^2}  & 1 & \cdot & \cdot & \cdot\\
        s_{21^4} & s_{1^2}  & 1 &  \cdot & \cdot \\
        s_{21^5} & s_{1^3} & s_1 &  1 & \cdot\\
        s_{21^7} & s_{1^5} & s_{1^3} & s_{1^2} & s_1
        \end{pmatrix}$\\
\hline
        $\begin{array}{c}
         \young{
        5 \cr
        4\cr
        3\cr
        2\cr
        1 \cr
        0 \cr
        {-1}\cr
        -2\cr
        -3\cr}\\
         \{[5,5],[3,4],[2,1],
         {[0,-1]},
         {[-1,-2]},[-3,-3]\}
         \end{array}$ &
        $\begin{pmatrix}
        s_1 & 1 & \cdot & \cdot & \cdot & \cdot\\
        s_{1^3}& s_{1^2}  & \cdot & \cdot & \cdot & \cdot \\
        s_{1^4} & s_{1^3} & 1 & \cdot & \cdot & \cdot\\
        s_{1^6} & s_{1^5} & s_{1^2}  & 1 &  \cdot & \cdot \\
        s_{1^7} & s_{1^6} & s_{1^3} & s_1 &  1 & \cdot\\
        s_{1^9} & s_{1^8} & s_{1^5} & s_{1^3} & s_{1^2} & s_1
        \end{pmatrix}$\\
        \hline
    \end{tabular}
    \end{center}
\end{table}


\begin{thebibliography}{9}
\bibitem{CYY} W. Y. C. Chen, G.-G. Yan, A. L. B. Yang,
Transformations of border strips and schur function determinants,
\textit{J. Algebraic Combin.}, to appear, arXiv:math.CO/0406250.



\bibitem{HG} A. M. Hamel, I. P. Goulden,
Planar decompositions of tableaux and Schur function determinants,
\textit{Europ. J. Combin.} \textbf{16} (1995), 461-477.

 \bibitem{LP2} A. Lascoux and P. Pragacz, Ribbon Schur functions,
 \textit{Europ. J. Combin.}, 9 (1988), 561-574.


\bibitem{KupJct} G. Kuperberg, Symmetries of plane partitions and
the permanent-determinant method, \textit{J. Combin. Theory Ser.
A} \textbf{68} (1994), no. 1, 115-151.

\bibitem{KupEjc} G. Kuperberg, Kasteleyn cokernels, \textit{Electron. J. Combin.} \textbf{9}
(2002), \#R29.


\bibitem{Sch} M. P. Sch\"utzenberger, La correspondance de Robinson,
in \textit{Combinatorics et repr\'esentation du groupe
symm\'etrique} (Table Ronde, Strasbourg 1976, D. Foata, ed.),
Lecture Notes in Math., Vol. 579, Springer, 1977, pp. 59-113.

\bibitem{Stanley} R. P. Stanley, Enumerative Combinatorics, Vol. 2,
Cambridge University Press, New York/Cambridge, 1999.

\bibitem{Z} A. Zelevinsky,
A generalization of the Littlewood-Richardson rule, \textit{J. Algebra}, 69 (1981), 82-94.

\end{thebibliography}
\end{document}